\documentclass[a4paper,12pt]{article}

\usepackage[all]{xy}
\usepackage[T1]{fontenc}
\usepackage[dvips]{graphicx}
\usepackage{amsmath,amssymb,amsfonts,amsthm,latexsym,mathrsfs,textcomp,verbatim, enumerate, soul, color}
\xyoption{web}
\title{Topologies for intermediate logics}
\author{Olivia Caramello \vspace{3 mm}\\ {\small DPMMS, University of Cambridge,}\\{\small Wilberforce Road, Cambridge CB3 0WB, U.K.}\\{\small O.Caramello@dpmms.cam.ac.uk}\thanks{The author gratefully acknowledges the support of a Research Fellowship from Jesus College, Cambridge (U.K.)}}

\date{May 11, 2012}
\begin{document}

%  change some single-character symbols to be more appropriate for logic
\mathcode`\<="4268  % < = \langle
\mathcode`\>="5269  % > = \rangle
\mathcode`\.="313A  % make . a binary  relation
\mathchardef\semicolon="603B % the original
\mathchardef\gt="313E
\mathchardef\lt="313C

\newcommand{\app}% application
 {{\sf app}}

\newcommand{\Ass}% category of assemblies
 {{\bf Ass}}

\newcommand{\ASS}% indexed version of \Ass
 {{\mathbb A}{\sf ss}}

\newcommand{\Bb}%blackboard bold
{\mathbb}

\newcommand{\biimp}% bi-implication
 {\!\Leftrightarrow\!}

\newcommand{\bim}% bimorphism
 {\rightarrowtail\kern-1em\twoheadrightarrow}

\newcommand{\bjg}% bi-judgement
 {\mathrel{{\dashv}\,{\vdash}}}

\newcommand{\bstp}[3]% bimorphism version of \stp
 {\mbox{$#1\! : #2 \bim #3$}}

\newcommand{\cat}%concatenation
 {\!\mbox{\t{\ }}}

\newcommand{\cinf}%C-infinity
 {C^{\infty}}

\newcommand{\cinfrg}%category of C-infinity rings
 {\cinf\hy{\bf Rng}}

\newcommand{\cocomma}[2]% cocomma object
 {\mbox{$(#1\!\uparrow\!#2)$}}

\newcommand{\cod}% codomain
 {{\rm cod}}

\newcommand{\comma}[2]% comma object
 {\mbox{$(#1\!\downarrow\!#2)$}}

\newcommand{\comp}% composition
 {\circ}

\newcommand{\cons}% concatenation
 {{\sf cons}}

\newcommand{\Cont}% category of continuous G-sets
 {{\bf Cont}}

\newcommand{\ContE}% continuous G-sets relative to $\cal E$
 {{\bf Cont}_{\cal E}}

\newcommand{\ContS}% ditto for $\cal S$
 {{\bf Cont}_{\cal S}}

\newcommand{\cover}% cover
 {-\!\!\triangleright\,}

\newcommand{\cstp}[3]% cover version of \stp
 {\mbox{$#1\! : #2 \cover #3$}}

\newcommand{\Dec}% decalage
 {{\rm Dec}}

\newcommand{\DEC}% decalage (\Bbb version)
 {{\mathbb D}{\sf ec}}

\newcommand{\den}[1]% denotation of #1 
 {[\![#1]\!]}

\newcommand{\Desc}% category of descent data
 {{\bf Desc}}

\newcommand{\dom}% domain
 {{\rm dom}}

\newcommand{\Eff}% effective topos
 {{\bf Eff}}

\newcommand{\EFF}% indexed version of \Eff
 {{\mathbb E}{\sf ff}}

\newcommand{\empstg}% empty string
 {[\,]}

\newcommand{\epi}% epimorphism
 {\twoheadrightarrow}

\newcommand{\estp}[3]% epimorphism version of \stp
 {\mbox{$#1 \! : #2 \epi #3$}}

\newcommand{\ev}% evaluation
 {{\rm ev}}

\newcommand{\Ext}% category of extracts
 {{\rm Ext}}

\newcommand{\fr}% Fraktur (i.e. Gothic)
 {\sf}

\newcommand{\fst}% first projection
 {{\sf fst}}

\newcommand{\fun}[2]% function-type
 {\mbox{$[#1\!\to\!#2]$}}

\newcommand{\funs}[2]% function-type as subscript
 {[#1\!\to\!#2]}

\newcommand{\Gl}% topos obtained by glueing
 {{\bf Gl}}

\newcommand{\hash}% hash sign (used as infix)
 {\,\#\,}

\newcommand{\hy}% hyphen (in math mode)
 {\mbox{-}}

\newcommand{\im}% image
 {{\rm im}}

\newcommand{\imp}% implication
 {\!\Rightarrow\!}

\newcommand{\Ind}[1]% ind-completion of #1
 {{\rm Ind}\hy #1}

\newcommand{\iten}[1]% enumerated item
{\item[{\rm (#1)}]}

\newcommand{\iter}% iterator
 {{\sf iter}}

\newcommand{\Kalg}%category of $K$-algebras
 {K\hy{\bf Alg}}

\newcommand{\llim}% left (inverse) limit
 {{\mbox{$\lower.95ex\hbox{{\rm lim}}$}\atop{\scriptstyle
{\leftarrow}}}{}}

\newcommand{\llimr}% right (inverse) limit
 {{\mbox{$\lower.95ex\hbox{{\rm lim}}$}\atop{\scriptstyle
{\rightarrow}}}{}}

\newcommand{\llimd}% \displaymath version of \llim
 {\lower0.37ex\hbox{$\pile{\lim \\ {\scriptstyle
\leftarrow}}$}{}}

\newcommand{\Mf}%category of manifolds
 {{\bf Mf}}

\newcommand{\Mod}% category of modest assemblies
 {{\bf Mod}}

\newcommand{\MOD}% indexed version of \Mod
{{\mathbb M}{\sf od}}

\newcommand{\mono}% monomorphism 
 {\rightarrowtail}

\newcommand{\mor}% class of morphisms
 {{\rm mor}}

\newcommand{\mstp}[3]% monomorphism version of \stp
 {\mbox{$#1\! : #2 \mono #3$}}

\newcommand{\Mu}%capital mu
 {{\rm M}}

\newcommand{\name}[1]% name of a relation
 {\mbox{$\ulcorner #1 \urcorner$}}

\newcommand{\names}[1]%\name used as subscript
 {\mbox{$\ulcorner$} #1 \mbox{$\urcorner$}}

\newcommand{\nml}% normal subgroup
 {\triangleleft}

\newcommand{\ob}% class of objects
 {{\rm ob}}

\newcommand{\op}% opposite category
 {^{\rm op}}
 
\newcommand{\palrr}[4]{ 
  \def\labelstyle{\scriptstyle} 
  \xymatrix{ {#1} \ar@<0.5ex>[r]^{#2} \ar@<-0.5ex>[r]_{#3} & {#4} } } 
  
\newcommand{\palrl}[4]{ 
  \def\labelstyle{\scriptstyle} 
  \xymatrix{ {#1} \ar@<0.5ex>[r]^{#2}  &  \ar@<0.5ex>[l]^{#3} {#4} } }  

\newcommand{\pepi}% partial epimorphism
 {\rightharpoondown\kern-0.9em\rightharpoondown}

\newcommand{\pmap}% partial map arrow
 {\rightharpoondown}

\newcommand{\Pos}% positivization of a coherent category
 {{\bf Pos}}

\newcommand{\prarr}% parallel pair of arrows
 {\rightrightarrows}

\newcommand{\princfil}[1]% principal filter
 {\mbox{$\uparrow\!(#1)$}}

\newcommand{\princid}[1]% principal ideal
 {\mbox{$\downarrow\!(#1)$}}

\newcommand{\prstp}[3]% parallel-pair version of \stp
 {\mbox{$#1\! : #2 \prarr #3$}}

\newcommand{\pstp}[3]% partial-map version of \stp
 {\mbox{$#1\! : #2 \pmap #3$}}

\newcommand{\relarr}% relation-type arrow
 {\looparrowright}

\newcommand{\rlim}% right limit, i.e. colimit
 {{\mbox{$\lower.95ex\hbox{{\rm lim}}$}\atop{\scriptstyle
{\rightarrow}}}{}}

\newcommand{\rlimd}% \displaymath version of \rlim
 {\lower0.37ex\hbox{$\pile{\lim \\ {\scriptstyle
\rightarrow}}$}{}}

\newcommand{\rstp}[3]% relation version of \stp
 {\mbox{$#1\! : #2 \relarr #3$}}

\newcommand{\scn}% Sierpinski cone
 {{\bf scn}}

\newcommand{\scnS}% ditto relative to $\cal S$
 {{\bf scn}_{\cal S}}

\newcommand{\semid}% semidirect product
 {\rtimes}

\newcommand{\Sep}% category of separated objects
 {{\bf Sep}}

\newcommand{\sep}% category of separated objects
 {{\bf sep}}

\newcommand{\Set}% category of sets
 {{\bf Set}}

\newcommand{\Sh}% category of sheaves
 {{\bf Sh}}

\newcommand{\ShE}% sheaves relative to $\cal E$
 {{\bf Sh}_{\cal E}}

\newcommand{\ShS}% ditto for $\cal S$
 {{\bf Sh}_{\cal S}}

\newcommand{\sh}% category of sheaves
 {{\bf sh}}

\newcommand{\Simp}% the simplicial category
 {{\bf \Delta}}

\newcommand{\snd}% second projection
 {{\sf snd}}

\newcommand{\stg}[1]% string of #1
 {\vec{#1}}

\newcommand{\stp}[3]% source-target predicate
 {\mbox{$#1\! : #2 \to #3$}}

\newcommand{\Sub}% subobject lattice
 {{\rm Sub}}

\newcommand{\SUB}% indexed category of subobjects
 {{\mathbb S}{\sf ub}}

\newcommand{\tbel}% totally below
 {\prec\!\prec}

\newcommand{\tic}[2]%term-in-context, etc.
 {\mbox{$#1\!.\!#2$}}

\newcommand{\tp}% is of type
 {\!:}

\newcommand{\tps}% subscript version of \tp
 {:}

\newcommand{\tsub}% truncated subtraction
 {\pile{\lower0.5ex\hbox{.} \\ -}}

\newcommand{\wavy}% wavy arrow
 {\leadsto}

\newcommand{\wavydown}% wavy downarrow
 {\,{\mbox{\raise.2ex\hbox{\hbox{$\wr$}
\kern-.73em{\lower.5ex\hbox{$\scriptstyle{\vee}$}}}}}\,}

\newcommand{\wbel}% way-below relation
 {\lt\!\lt}

\newcommand{\wstp}[3]% wavy version of \stp
 {\mbox{$#1\!: #2 \wavy #3$}}
 
\newcommand{\fu}[2]
{[#1,#2]}

%\newcommand{\st}[2]% source-target predicate
% {\mbox{$#1 \to #2$}} 

%===========================================================================
%	END OF PROOF BOX
%
%
%  The complexity of the macro necessary to get a little box on the
%  right-hand-side at the end of a proof is amazing.  It really does
%  have to be this long!  Otherwise you're liable to get it at the
%  beginning of the next line, or even on the next page.
%
\def\pushright#1{{%        set up
   \parfillskip=0pt            % so \par doesnt push \square to left
   \widowpenalty=10000         % so we dont break the page before \square
   \displaywidowpenalty=10000  % ditto
   \finalhyphendemerits=0      % TeXbook exercise 14.32
  %
  %                 horizontal
   \leavevmode                 % \nobreak means lines not pages
   \unskip                     % remove previous space or glue
   \nobreak                    % don't break lines
   \hfil                       % ragged right if we spill over
   \penalty50                  % discouragement to do so
   \hskip.2em                  % ensure some space
   \null                       % anchor following \hfill
   \hfill                      % push \square to right
   {#1}                        % the end-of-proof mark (or whatever)
  %
  %                   vertical
   \par}}                      % build paragraph

% prefer proofs with statements, also space after
\def\qed{\pushright{$\square$}\penalty-700 \smallskip}

\newtheorem{theorem}{Theorem}[section]

\newtheorem{proposition}[theorem]{Proposition}

\newtheorem{scholium}[theorem]{Scholium}

\newtheorem{lemma}[theorem]{Lemma}

\newtheorem{corollary}[theorem]{Corollary}

\newtheorem{conjecture}[theorem]{Conjecture}

\newenvironment{proofs}%
 {\begin{trivlist}\item[]{\bf Proof }}%
 {\qed\end{trivlist}}

  \newtheorem{rmk}[theorem]{Remark}
\newenvironment{remark}{\begin{rmk}\em}{\end{rmk}}

  \newtheorem{rmks}[theorem]{Remarks}
\newenvironment{remarks}{\begin{rmks}\em}{\end{rmks}}

  \newtheorem{defn}[theorem]{Definition}
\newenvironment{definition}{\begin{defn}\em}{\end{defn}}

  \newtheorem{eg}[theorem]{Example}
\newenvironment{example}{\begin{eg}\em}{\end{eg}}

  \newtheorem{egs}[theorem]{Examples}
\newenvironment{examples}{\begin{egs}\em}{\end{egs}}

%%%%%%%%%%%%%%%%%%%%%%%%%%%%%%%%%%%%%%%%%%%%%%%%%%%%%%%%%%%%%%%%%%%%%%

\bgroup           % fake a titlepage 
\let\footnoterule\relax  % no rule above thanks footnotes 
\maketitle

\begin{abstract}
We investigate the problem of characterizing the classes of Grothen-\\dieck toposes whose internal logic satisfies a given assertion in the theory of Heyting algebras, and introduce natural analogues of the double negation and De Morgan topologies on an elementary topos for a wide class of intermediate logics.  
\end{abstract} 
\egroup

\tableofcontents

\newpage
\section{Introduction}

In light of the fact that the internal logic of a topos is at least intuitionistic, it is natural to investigate the class of toposes whose logic satisfies some additional assertion written in the theory of Heyting algebras (in the sense that such assertion is satisfied in the internal Heyting algebra to the topos given by its subobject classifier). In fact, besides its natural theoretical interest, such an investigation can pave the way for the introduction of new topos-theoretic invariants admitting bijective site characterizations, something which, among the other things, is particularly relevant in connection to the methodology `toposes as bridges' of \cite{OC10}. 

A related problem is the construction of a universal way of associating to a general topos a subtopos of it satisfying a given intermediate logic. For example, the subtopos of a given elementary topos consisting of its double-negation sheaves can be seen as a universal way of making the topos Boolean, as it can be characterized as the largest dense Boolean subtopos of the given topos; similarly, the subtopos of a given elementary topos consisting of its sheaves with respect to the De Morgan topology (as introduced in \cite{OC3}) can be characterized as its largest dense subtopos satisfying De Morgan's law. These concepts have proved to be fruitful in different contexts (cf. for example \cite{OC2} and \cite{OC4}), so it is natural to look for analogues of them for general intermediate logics.

To this end, we identify a stronger property enjoyed by the Booleanization (resp. DeMorganization) of a topos, namely the fact that these subtoposes are not only the largest among the dense subtoposes satisfying the law of excluded middle (resp. De Morgan's law), as shown in \cite{OC3}, but more generally among all the subtoposes with the property that their associated sheaf functor preserves the pseudocomplementation operation on subobjects. This remark indicates that, for any intermediate logic whose definition involves, besides the conjunction and disjunction connectives, the connective $\neg$ (resp. the connective $\imp$), it is natural to look for a dense subtopos of a given topos satisfying that logic and containing all the subtoposes of the given topos which satisfy that logic and whose corresponding associated sheaf functor preserves the pseudocomplementation (resp. Heyting implication) operation on subobjects. In fact, we prove that for a wide class of intermediate logics such a subtopos exists, and provide an explicit description of it in the case of a general topos of sheaves on a site, as well as in the particular cases of localic and presheaf toposes; in particular, since subtoposes of localic toposes are localic, this construction also yields, for any given locale $L$, a dense sublocale of it satisfying the intermediate logic in question and containing any sublocale of $L$ which satisfies the logic and whose quotient map preserves the pseudocomplementation (resp. the Heyting implication) operation.    

The plan of the paper is as follows.

In section \ref{sec2} we discuss the problem of finding explicit criteria for (the subobject classifier of) a Grothendieck topos to satisfy a given first-order sequent written in the theory of Heyting algebras; in particular, this leads to criteria for a Grothendieck topos to satisfy a given intermediate logic. We specifically address the case of presheaf toposes, localic toposes and classifying toposes (i.e., toposes of sheaves on the syntactic site of a geometric theory), establishing appropriate criteria for them.  

In section \ref{sec3} we study the local operators on elementary toposes with the property that their corresponding associated sheaf functor preserves the pseudocomplementation (resp. the implication) operation on subobjects. This paves the way for the introduction, carried out in section \ref{sec4}, of appropriate analogues of the double-negation and De Morgan topologies on an elementary topos for a wide class of intermediate logics.

In section \ref{sec4}, besides introducing these new constructions, we discuss the problem of calculating them in several cases of interest.

\section{Criteria for a Grothendieck topos to satisfy an intermediate logic}\label{sec2}

For any first-order sequent $\sigma$ written in the theory of Heyting algebras, it is possible to obtain explicit criteria for the subobject classifier $\Omega_{\Sh({\cal C}, J)}$ of a topos $\Sh({\cal C}, J)$ of sheaves on a site $({\cal C}, J)$ to satisfy $\sigma$, by using the following explicit (and easily provable) descriptions of the internal Heyting operations on $\Omega_{\Sh({\cal C}, J)}$, together with the equally explicit descriptions of the interpretation of first-order connectives and quantifiers in a Grothendieck topos, as given in section III.8 of \cite{MM}:

\begin{itemize}
\item the subobject classifier $\Omega_{\Sh({\cal C}, J)}$ of $\Sh({\cal C}, J)$ is defined by the following formulas:
\[
\Omega_{\Sh({\cal C}, J)}(c)=\textrm{\{R | R is a $J$-closed sieve on $c$\}}
\]
for any object $c\in \cal C$,
\[
\Omega_{\Sh({\cal C}, J)}(f)=f^{\ast}(-)
\]
for any arrow $f$ in $\cal C$, where $f^{\ast}$ denotes the operation of pullback of sieves in $\cal C$ along $f$;

\item the bottom element
\[
0:1\to \Omega_{\Sh({\cal C}, J)}
\]
of the algebra $\Omega$ is defined by setting $0(c)(\ast)$ equal to the $J$-closure 
\[
\overline{\emptyset_{c}}^{J}=\{f:d\to c \textrm{ | } \emptyset \in J(d)\}
\]
of the empty sieve on $c$ (for any object $c$ of $\cal C$);
\item the top element
\[
1:1\to \Omega_{\Sh({\cal C}, J)}
\]
of the algebra $\Omega_{\Sh({\cal C}, J)}$ is defined by setting $1(c)(\ast)$ equal to the maximal sieve on $c$;
\item the meet operation 
\[
\wedge: \Omega_{\Sh({\cal C}, J)}\times \Omega_{\Sh({\cal C}, J)} \to \Omega_{\Sh({\cal C}, J)}
\]
on $\Omega_{\Sh({\cal C}, J)}$ is given by the formula
\[
\wedge(c)(S, T)=S\cap T
\] 
(for any object $c$ of ${\cal C}$ and any $J$-closed sieves $S$ and $T$ on $c$);
\item the join operation
\[
\vee: \Omega_{\Sh({\cal C}, J)}\times \Omega_{\Sh({\cal C}, J)} \to \Omega_{\Sh({\cal C}, J)}
\]
on $\Omega_{\Sh({\cal C}, J)}$ is given by the formula 
\[
\vee(c)(S, T)=\{f:d\to c \textrm{ in ${\cal C}$ | } f^{\ast}(S\cup T)\in J(d)\}
\]
(for any object $c$ of ${\cal C}$ and any $J$-closed sieves $S$ and $T$ on $c$);
\item the Heyting implication operation
\[
\imp: \Omega_{\Sh({\cal C}, J)}\times \Omega_{\Sh({\cal C}, J)} \to \Omega_{\Sh({\cal C}, J)}
\]
is defined by the formula
\[
\imp(c)(S, T)=\{f:d\to c \textrm{ in ${\cal C}$ | } f^{\ast}(S)\subseteq f^{\ast}(T)\} 
\]
(for any object $c$ of ${\cal C}$ and any $J$-closed sieves $S$ and $T$ on $c$); 
\item the Heyting pseudocomplementation operation
\[
\neg:\Omega_{\Sh({\cal C}, J)} \to \Omega_{\Sh({\cal C}, J)}
\] 
is given by the formula
\[
\neg(c)(S)=\{f:d\to c \textrm{ | for all } g:e\to d, \textrm{  $f\circ g\in S$ implies } \emptyset \in J(e)\}
\] 
\end{itemize}

A mechanical application of these formulas allows us to achieve completely explicit criteria for $\Omega_{\Sh({\cal C}, J)}$ to satisfy any first-order sequent $\sigma$ in the theory of Heyting algebras, of the form `$\sigma$ is satisfied in $\Omega_{\Sh({\cal C}, J)}$ if and only if the site $({\cal C}, J)$ satisfies a property $P_{({\cal C}, J)}$ explicitly written in the language of the site $({\cal C}, J)$'. Notice that, in light of the fact that the subobject classifier of a topos is a topos-theoretic invariant, such criteria can be profitably applied in connection to the philosophy `toposes as bridges' of \cite{OC10}. 

Before proceeding to a selection of examples of such criteria, let us establish some general results enabling us to obtain, in a variety of naturally occurring situations, simplifications of them.  
  
The following result provides a relationship between the notion of validity of a cartesian sequent in the theory of Heyting algebras in the internal Heyting algebra of a topos given by its subobject classifier and the concept of validity of the sequent in the `external' subobject lattices in the topos.  

\begin{theorem}\label{cartesian}
Let $\cal E$ be a locally small elementary topos, $\sigma$ be a cartesian (in particular, Horn) sequent in the theory of Heyting algebras, and $\cal C$ be a set of objects of $\cal E$ such that the class of objects of $\cal E$ which admit a monomorphism to an object of $\cal C$ form a separating set for $\cal E$. Then $\sigma$ is valid in the internal algebra $\Omega_{{\cal E}}$ of the topos $\cal E$ given by its subobject classifier if and only if it is valid in $\Sub_{\cal E}(c)$ for every $c\in {\cal C}$ (where this poset is regarded as a model of the theory of Heyting algebras).
\end{theorem}

\begin{proofs}
For any locally small topos $\cal E$, the Yoneda embedding $y:{\cal E}\to [{\cal E}^{\textrm{op}}, \Set]$ is a cartesian functor, whence $y$ preserves and the interpretation of all the cartesian formulae. From this it follows that, given any cartesian sequent $\sigma$ in the theory of Heyting algebras, the internal Heyting algebra $\Omega_{{\cal E}}$ in $\cal E$ given by its subobject classifier satisfies $\sigma$ if and only if every frame $\Sub_{\cal E}(e)\cong Hom_{\cal E}(e, \Omega)$ in $\cal E$ satisfies $\sigma$. Now, given an object $e\in \cal E$, if ${\cal C}'$ is a separating set for $\cal E$ then $e$ can be expressed as a quotient of a coproduct of objects in ${\cal C}'$, that is there exists a set-indexed family $\{c_{i} \textrm{ | } i\in I\}$ of objects in ${\cal C}'$ and an epimorphism $p:\coprod_{i\in I}c_{i}\epi e$; so the pullback functor $p^{\ast}:\Sub_{\cal E}(e) \rightarrow \Sub_{\cal E}(\coprod_{i\in I}c_{i})\cong \prod_{i\in I}\Sub_{\cal E}(c_{i})$ is logical and conservative (cf. Example A4.2.7(a) \cite{El}) and hence $\Sub_{\cal E}(e)$ satisfies a first-order sequent $\sigma$ if all the $\Sub_{\cal E}(c_{i})$ do. On the other hand, if $m:b\mono a$ is a monomorphism in $\cal E$ then the pullback functor $m^{\ast}:\Sub_{\cal E}(a)\to \Sub_{\cal E}(b)$ is logical and essentially surjective (since for any subobject $k:c\mono b$, $k\cong m^{\ast}(m\circ k)$); so, if $\Sub_{\cal E}(a)$ satisfies $\sigma$ then $\Sub_{\cal E}(b)$ satisfies $\sigma$. Our thesis now follows immediately from the combination of these two facts.    
\end{proofs}

This result has a couple of useful corollaries. 

\begin{corollary}
Let $\mathbb T$ be a geometric theory over a signature $\Sigma$. Then the classifying topos of $\mathbb T$ internally satisfies a cartesian sequent $\sigma$ in the theory of Heyting algebras if and only if all the frames of $\mathbb T$-provable equivalence classes of geometric formulae over $\Sigma$ in a given context satisfy $\sigma$.
\end{corollary}

\begin{proofs}
It suffices to observe that, via the Yoneda embedding $y:{\cal C}_{\mathbb T}\to \Sh({\cal C}_{\mathbb T}, J_{\mathbb T})$, the set of objects $\cal C$ of the form $y(\{\vec{x}. \top\})$ (for any context $\vec{x}$) satisfies the hypotheses of the theorem, and the subobjects in $\Sh({\cal C}_{\mathbb T}, J_{\mathbb T})$ of an object of the form $y(\{\vec{x}. \top\})$ can be identified with the $\mathbb T$-provable equivalence classes of geometric formulae over $\Sigma$ in the context $\vec{x}$ (by Lemma D1.4.4(iv) \cite{El}). 
\end{proofs}

Note that this corollary can be applied in particular in the context of the investigation of the Lee identities on the classifying topos of a geometric theory.

The following result shows that, for any cartesian sequent in the theory of Heyting algebras, its internal validity in a localic topos $\Sh(L)$ is equivalent to its `external' validity in the locale $L$.  

\begin{corollary}
Let $L$ be a locale and $\sigma$ be a cartesian (in particular, Horn) sequent in the theory of Heyting algebras. Then $\sigma$ is valid in the algebra $\Omega_{\Sh(L)}$ of the topos $L$ if and only if it is valid in $L$ (considered as a model of the theory of Heyting algebras). 
\end{corollary}

\begin{proofs}
This follows immediately from Theorem \ref{cartesian} observing that the family $\cal C$ consisting of the terminal object of $\Sh(L)$ satisfies its hypotheses.
\end{proofs}

\begin{remark}
These corollaries notably apply to all the intermediate logics; indeed, every intermediate logic can be seen as a Horn theory which extends the theory of Heyting algebras over its signature, whose axioms are all of the form $(\top \vdash_{\vec{x}} t_{1}=t_{2})$ where $t_{1}$ and $t_{2}$ are two terms in the context $\vec{x}$ written in the language of Heyting algebras.
\end{remark}

Let us conclude this section with some examples of criteria for a Grothen-\\-dieck topos to satisfy an intermediate logic, obtained through an application of the general method described above. In \cite{OC3}, site characterizations for the property of a Grothendieck topos to be Boolean (resp. De Morgan), were obtained. Another interesting example is given by G\"{o}del-Dummett logic, that is the logic obtained from intuitionistic propositional logic by adding the axiom scheme $(p\imp q) \vee (q\imp p)$. It is easy to calculate, by applying the formulas established above, that a Grothendieck topos $\Sh({\cal C}, J)$ satisfies G\"{o}del-Dummett logic if and only if for every $J$-closed sieves $R$ and $S$ on an object $c\in {\cal C}$, the sieve $\{f:d\to c \textrm{ | } f^{\ast}(R)\subseteq f^{\ast}(S) \textrm{ or } f^{\ast}(S)\subseteq f^{\ast}(R)\}$ is $J$-covering. 

In the case of presheaf toposes, these criteria for a Grothendieck topos to be Boolean (resp. to be De Morgan, to satisfy G\"{o}del-Dummett logic) specialize to the following well-known results (cf. \cite{El} for the first two and \cite{morgan} for the third):   

\begin{proposition}
Let $\cal C$ be a small category. Then 

\begin{enumerate}[(i)]
\item the topos $[{\cal C}^{\textrm{op}}, \Set]$ is Boolean if and only if $\cal C$ is a groupoid;

\item the topos $[{\cal C}^{\textrm{op}}, \Set]$ is De Morgan if and only if $\cal C$ satisfies the right Ore condition (i.e., for any two arrows $f:b\to a$ and $g:c\to a$ with common codomain there is an object $d$ and arrows $h:d\to b$ and $k:d\to c$ such that $f\circ h=g\circ k$);

\item the topos $[{\cal C}^{\textrm{op}}, \Set]$ satisfies G\"{o}del-Dummett logic if and only if $\cal C$ satisfies the following property: for any arrows $f:b\to a$ and $g:c\to a$ with common codomain, either $f$ factors through $g$ or $g$ factors through $f$.
\end{enumerate}
\end{proposition}

This proposition shows that these invariants, when considered on a pre-\\-sheaf topos $[{\cal C}^{\textrm{op}}, \Set]$, capture interesting geometrical properties of the category $\cal C$; on the other hand, considered on a topos $\Sh(X)$ of sheaves on a topological space $X$, they specialize to important properties of $X$ (namely, the property of $X$ to be almost discrete, resp. extremally disconnected, satisfying the property that the closure of any open set is an extremally disconnected, cf. \cite{El} and \cite{morgan}). This accounts for the unifying power of these invariants (and more generally of those given by the interpretation of first-order sequents in the theory of Heyting algebras in the subobject classifier of the topos), in the sense that they can be effectively used in presence of any Morita-equivalence of toposes
\[
\Sh({\cal C}, J)\simeq \Sh({\cal D}, K)
\]
to operate an \emph{automatic} transfer of properties between the sites $({\cal C}, J)$ and $({\cal D}, K)$ (cf. \cite{OC11} for a selection of applications of this general method). 

As another example, let us consider the intermediate logic known as Kreisel-Putnam logic, that is the logic obtained from intuitionistic propositional logic by adding the axiom scheme $(\neg p \imp (q\vee r)) \imp ((\neg p \imp q) \vee (\neg p \imp r))$.

It is immediate to see, by using the formulas established above, that a Grothendieck topos $\Sh({\cal C}, J)$ satisfies this logic if and only if for any $J$-closed sieves $R$, $S$ and $T$ on an object $c$ of $\cal C$, if $\neg R=\{f:d\to c \textrm{ | } f^{\ast}(R)=\emptyset \}$ is equal to $\{f:d\to c \textrm{ | } f^{\ast}(S)\in J(d) \textrm{ or } f^{\ast}(T)\in J(d)\}$ then $\neg R=S$ or $\neg R=T$. In particular, a presheaf topos $[{\cal C}^{\textrm{op}}, \Set]$ satisfies this logic if and only if every stably non empty sieve (that is, any sieve $R$ on an object $c$ such that for any arrow $f$ with codomain $c$, $f^{\ast}(S)\neq \emptyset$) is \emph{indecomposable}, in the sense that for any two sieves $S$ and $T$ on $\cal C$, $R=S\cup T$ implies either $R=S$ or $R=T$. 

The following result provides a characterization of indecomposable sieves.

\begin{proposition}\label{indecomp}
Let $\cal C$ be a category and $R$ be a sieve in $\cal C$ on an object $c$. Then $R$ is indecomposable if and only if it satisfies the following property: for any arrows $f:d\to c$ and $g:e\to c$ in $R$ there exists $h:a\to c$ in $R$ such that both $f$ and $g$ factor through $h$.
\end{proposition}

\begin{proofs}
Let us suppose that $R$ is indecomposable. Given $f:d\to c$ in $R$, define $T_{f}:=\{g:e\to c \in R  \textrm{ | $f$ does not factor through $g$}\}$ and $H_{f}:=\{g:e\to c\in R  \textrm{ | $f$ factors through $g$}\}$. Clearly, $T_{f}$ is a sieve; if $\overline{H_{f}}$ is the sieve generated by $H_{f}$ then we have the decomposition $R=T_{f}\cup \overline{H_{f}}$ of $R$ as a union of two sieves on $c$. Since $f\notin T_{f}$, from the fact that $R$ is indecomposable it follows that $R=\overline{H_{f}}$. So, for any arrow $g:e\to c$ in $R$, $g$ belongs to $\overline{H_{f}}$, that is there exists $h\in H_{f}$ through which both $g$ and $f$ factor.

Conversely, suppose that $R$ satisfies the property stated in the Proposition; we want to prove that $R$ is indecomposable. Let $R=S\cup T$. Suppose for contradiction that $R\neq S$ and $R\neq T$; then in particular $S$ is not a subset of $T$ and $T$ is not a subset of $S$. So there exists $f\in T$ such that $f$ does not belong to $S$ and $g\in S$ such that $g$ does not belong to $T$. Now, by our hypothesis there exists $h\in R$ such that both $f$ and $g$ factor through $h$. Since $R=S\cup T$ then $h$ must belong either to $S$ or to $T$, leading in either case to a contradiction.   
\end{proofs}

\section{Dense, weakly open and implicationally open subtoposes}\label{sec3}

In this section we investigate the local operators $j$ on an elementary topos $\cal E$ such that the associated sheaf functor $a_{j}:{\cal E}\to \sh_{j}({\cal E})$ preserves the pseudocomplementation operation (resp. the Heyting implication operation) on subobjects. 

We denote by $c_{j}$ the closure operation on subobjects of $\cal E$ corresponding to a local operator $j$ on $\cal E$, and we write $i_{j}:\sh_{j}({\cal E})\hookrightarrow {\cal E}$ for the obvious inclusion. We denote by $\Omega$ the (codomain of the) subobject classifier of $\cal E$, and by $\Omega_{j}$ the subobject classifier of $\sh_{j}({\cal E})$, given by the equalizer $e_{j}:\Omega_{j}\mono \Omega$ of the two arrows $j, 1_{\Omega}:\Omega\to \Omega$.

The following proposition gives a bunch of alternatives ways of characterizing the dense local operators on an elementary topos $\cal E$, i.e. the local operators $j$ on $\cal E$ such that $j\leq \neg\neg$. This result (in which the first four characterizations are well-known) is useful, among the other things, for illuminating the subtle relationship between dense operators and local operators whose associated sheaf functors preserves the pseudocomplementation operation on subobjects. 

\begin{proposition}
Let $\cal E$ be an elementary topos and $j$ be a local operator on $\cal E$. Then the following conditions are equivalent.
\begin{enumerate}[(i)]

\item $j\leq \neg\neg$;

\item $c_{j}(0\mono 1)=0\mono 1$;

\item For any monomorphism $m$ in $\cal E$, $a_{j}(m)\cong 0$ in $\sh_{j}({\cal E})$ implies $m\cong 0$ in $\cal E$;

\item The inclusion $i_{j}:\sh_{j}({\cal E})\hookrightarrow {\cal E}$ preserves the initial object.

\item The diagram
\[  
\xymatrix {
\Omega \ar[d]^{j} \ar[r]^{\neg} & \Omega \ar[d]^{j}\\
\Omega \ar[r]_{\neg} & \Omega}
\]
commutes;

\item either (equivalently, both) of the triangles 

\[  
\xymatrix {
\Omega \ar[d]^{j} \ar[dr]^{\neg} & &  \Omega \ar[dr]_{\neg} \ar[r]^{\neg}  & \Omega \ar[d]^{j} \\
\Omega \ar[r]_{\neg} & \Omega & & \Omega}
\]
commutes;

\item The closure operation $c_{j}$ preserves the pseudocomplementation operation on subobjects in $\cal E$;

\item The diagram
\[  
\xymatrix {
\Omega_{j} \ar[r]^{\neg_{j}} \ar[d]_{e_{j}} &  \Omega_{j} \ar[d]^{e_{j}} \\
\Omega \ar[r]_{\neg} & \Omega}
\]
commutes, where $\neg_{j}:\Omega_{j}\to \Omega_{j}$ is the pseudocomplementation operation in the internal Heyting algebra in $\sh_{j}({\cal E})$ given by its subobject classifier $\Omega_{j}$;

\item The equality of subtoposes
\[
\sh_{\neg\neg}(\sh_{j}({\cal E}))=\sh_{\neg\neg}({\cal E})
\]
holds.

\end{enumerate}
\end{proposition}

\begin{proofs}
The equivalence of the first four conditions is well-known. Condition $(vii)$ is clearly the `externalization' of condition $(v)$ and, as such, it is equivalent to it. Let us prove that $(i)$ implies $(vii)$. For any subobject $m$, if $j$ satisfies $(i)$ then $c_{j}(\neg m)\leq \neg c_{j}(m)$, while the converse inequality follows from the fact that $\neg c_{j}(m)=c_{j}(\neg c_{j}(m))\leq c_{j}(\neg m)$ (since $\neg c_{j}(m)$ is $\neg\neg$-closed and hence $c_{j}$-closed). Conversely, if $j$ satisfies $(vii)$ then $c_{j}(0\mono 1)=c_{j}(\neg 1_{1})=\neg c_{j}(1_{1})=0\mono 1$ and hence $(ii)$, equivalently $(i)$, is satisfied.

Let us now prove the equivalence of $(vi)$ and $(vii)$. Under the hypothesis that $(vii)$ (equivalently $(i)$) holds, the first triangle commutes since for any subobject $m$, $c_{j}(\neg m)=\neg m$ (since $\neg m$ is $\neg\neg$-closed and hence $c_{j}$-closed); the fact that the second triangle commutes follows immediately from the commutativity of the square in $(v)$ and of the first triangle. Conversely, if both triangles commute then the square in $(v)$ commutes, equivalently $(vii)$ holds.

To prove that $(ii)$ implies $(viii)$ we observe that if $0\mono 1$ is $c_{j}$-closed then its classifying map $\bot:1\to \Omega$ factors through $e_{j}:\Omega_{j}\mono \Omega$, and its factorization $\bot_{j}:1\to \Omega_{j}$ can be identified with the bottom element of the internal Heyting algebra in $\sh_{j}({\cal E})$ given by its subobject classifier $\Omega_{j}$; since for any local operator $j$ we have a commutative diagram
\[  
\xymatrix {
\Omega_{j}\times\Omega_{j} \ar[d]^{e_{j}\times e_{j}} \ar[rr]^{\imp_{j}} & & \Omega_{j} \ar[d]^{e_{j}} \\
\Omega\times\Omega \ar[rr]_{\imp} & & \Omega}
\] 
(cf. the proof of Proposition 6.8 \cite{OC6}), the square
\[  
\xymatrix {
\Omega_{j} \ar[r]^{\neg_{j}} \ar[d]_{e_{j}} &  \Omega_{j} \ar[d]^{e_{j}} \\
\Omega \ar[r]_{\neg} & \Omega}
\]
commutes. 

Conversely, let us suppose that the square in condition $(viii)$ commutes. If we denote by $\top_{j}$ the factorization of the true arrow $\top:1\to \Omega$ across the arrow $e_{j}:\Omega_{j}\to \Omega$ then the commutativity of this square implies that the classifying map $\bot=\neg\circ \top=\neg\circ e_{j}\circ \top_{j}:1\to \Omega$ of the subobject $0\mono 1$ factors through $e_{j}$, that is the subobject $0\mono 1$ is $c_{j}$-closed, as required.  
  
It remains to prove that $(ix)$ is equivalent to $(i)$. It is well-known that if $j$ satisfies $(i)$ then $\sh_{\neg\neg}(\sh_{j}({\cal E}))=\sh_{\neg\neg}({\cal E})$ (cf. for example the proof of Lemma A4.5.21 \cite{El}); conversely, if $\sh_{\neg\neg}(\sh_{j}({\cal E}))=\sh_{\neg\neg}({\cal E})$ then clearly $\sh_{\neg\neg}({\cal E})\subseteq \sh_{j}({\cal E})$ and hence $j\leq \neg\neg$.    
\end{proofs}

\begin{definition}
Let $\cal E$ be a topos and $j$ a local operator on $\cal E$, with associated sheaf functor $a_{j}:{\cal E}\to \sh_{j}({\cal E})$.

\begin{enumerate}[(a)]

\item We say that $j$ (resp. the subtopos $\sh_{j}({\cal E})$) is \emph{weakly open} if $a_{j}$ preserves the pseudocomplementation of subobjects; 

\item We say that $j$ (resp. the subtopos $\sh_{j}({\cal E})$) is \emph{implicationally open} if $a_{j}$ preserves the Heyting implication of subobjects.

\end{enumerate}
\end{definition}

Notice that every implicationally open local operator is weakly open (indeed, the associated sheaf functor always preserves the initial object, and for any object $a$ of $\cal E$, $\neg a=a\imp 0$).

Recall from \cite{OC7} (Proposition 6.3) that every dense local operator is weakly open. On the other hand, the converse does not hold, since every open local operator is weakly open but not in general dense. Indeed, if $o(U)$ is the open local operator on a topos $\cal E$ corresponding to a subterminal $U$ in $\cal E$, the corresponding closure operation $c_{o(U)}$ sends a subobject $A'\mono A$ to the implication $(A\times U)\imp A'$ in the Heyting algebra $\Sub_{\cal E}(A)$ (cf. \cite{El}) and hence $c_{o(U)}(0\mono 1)=\neg U$, which is in general different from $0$.

\begin{proposition}\label{dense}
Let $\cal E$ be an elementary topos and $j$ be a local operator on $\cal E$. Then the following conditions are equivalent.

\begin{enumerate}[(i)]
\item The associated sheaf functor $a_{j}:{\cal E}\to \sh_{j}({\cal E})$ preserves the Heyting implication of subobjects;

\item The diagram
\[  
\xymatrix {
\Omega\times\Omega \ar[d]^{j\times j} \ar[rr]^{\imp} & & \Omega \ar[d]^{j}\\
\Omega\times\Omega \ar[rr]_{\imp} & & \Omega}
\]
commutes.
\end{enumerate}
\end{proposition}

\begin{proofs}
We denote the Heyting implication of subobjects in $\cal E$ (resp. in $\sh_{j}({\cal E})$) by $\imp$ (resp. by $\imp_{j}$).

Recall from \cite{OC6} (cf. the proof of Proposition 6.8) that we have a commutative diagram  
\[  
\xymatrix {
\Omega_{j}\times\Omega_{j} \ar[d]^{e_{j}\times e_{j}} \ar[rr]^{\imp_{j}} & & \Omega_{j} \ar[d]^{e_{j}} \\
\Omega\times\Omega \ar[rr]_{\imp} & & \Omega}
\]
in $\cal E$, where $e_{j}:\Omega_{j}\mono \Omega$ is the equalizer of the pair of maps $j, 1_{\Omega}:\Omega\to \Omega$. This means that, given any two subobjects $m,n$ of a given object in $\sh_{j}({\cal E})$, their Heyting implication $m\imp_{j} n$ in $\sh_{j}({\cal E})$ coincides with their Heyting implication $m\imp n$ in $\cal E$.

Recall that for any subobject $m:A'\mono A$ in $\cal E$, we have a pullback diagram
 \[  
\xymatrix {
c_{j}(a') \ar[r]^{c_{j}(m)} \ar[d] & a \ar[d]^{\eta_{A}} \\
a_{j}(a') \ar[r]^{a_{j}(m)} & a_{j}(a) }
\]   
in $\cal E$, where $\eta_{A}$ is the unit of the reflection $\sh_{j}({\cal E}) \hookrightarrow {\cal E}$; in other words, $c_{j}(m)=\eta_{a}^{\ast}(a_{j}(m))$. Let us suppose that $(i)$ holds and derive $(ii)$. Notice that condition $(ii)$ is equivalent to saying that the closure operation $c_{j}$ preserves the Heyting implication of subobjects. Now, if condition $(i)$ holds then for any subobjects $m$ and $n$ of an object $a$ in $\cal E$, $c_{j}(m\imp n)=\eta_{a}^{\ast}(a_{j}(m\imp n))=\eta_{a}^{\ast}(a_{j}(m)\imp_{j} a_{j}(n))$ and hence, by the remark above, $c_{j}(m\imp n)=\eta_{a}^{\ast}(a_{j}(m)\imp a_{j}(n))=\eta_{a}^{\ast}(a_{j}(m)) \imp \eta_{a}^{\ast}(a_{j}(m))=c_{j}(m)\imp c_{j}(n)$. 

Conversely, let us prove that $(ii)$ implies $(i)$. If $m$ and $n$ are two subobjects of an object $a$ in $\cal E$, $a_{j}(m\imp n)=a_{j}(c_{j}(m\imp n))=a_{j}(c_{j}(m) \imp c_{j}(n))=a_{j}(\eta_{a}^{\ast}(a_{j}(m))\imp \eta_{a}^{\ast}(a_{j}(n)))=a_{j}(\eta_{a}^{\ast}(a_{j}(m) \imp a_{j}(n)))=a_{j}(\eta_{a}^{\ast}(a_{j}(m) \imp_{j} a_{j}(n)))\\$ $=a_{j}(m) \imp_{j} a_{j}(n)$, where the last passage follows from the fact that for any subobject $r:b\mono a_{j}(a)$ in $\sh_{j}({\cal E})$, $a_{j}(\eta_{a}^{\ast}(r))=r$.
\end{proofs}

The proposition shows that the class of local operators that we call implicationally open coincides with the class of operators characterized by Proposition A4.5.8 \cite{El}. We shall now establish a criterion for identifying implicationally open local operators on Grothendieck toposes. To prove it, we need a lemma.

\begin{lemma}
Let $\cal E$ be an elementary topos, $j$ be a local operation on it and $s:{\cal F}\hookrightarrow \sh_{j}({\cal E})$ be a subtopos of $\sh_{j}({\cal E})$. Then $s$ is implicationally open if and only if the diagram
\[  
\xymatrix {
\Omega_{j}\times\Omega_{j} \ar[d]^{(s\times s)\circ (e_{j}\times e_{j})} \ar[rr]^{\imp_{j}} & & \Omega_{j} \ar[d]^{s\circ e_{j}} \\
\Omega\times\Omega \ar[rr]_{\imp} & & \Omega}
\]
commutes, where $k$ is the local operator on $\cal E$ corresponding to the geometric inclusion given by composite of $s$ with the canonical inclusion $\sh_{j}({\cal E})\hookrightarrow {\cal E}$.    
\end{lemma}

\begin{proofs}
If $k'$ is the local operator on $\sh_{j}({\cal E})$ corresponding to the subtopos $i$ then the diagram
 \[  
\xymatrix {
\Omega_{j} \ar[d]^{e_{j}} \ar[r]^{k'} & \Omega_{j} \ar[d]^{e_{j}} \\
\Omega \ar[r]_{k}  & \Omega}
\]
commutes (by Proposition 6.4 \cite{OC6}). The thesis thus follows from the fact that $e_{j}$ is a monomorphism.
\end{proofs}

\begin{proposition}
Let $\Sh({\cal C}, J)$ be a Grothendieck topos, and $k$ be a local operator on $\Sh({\cal C}, J)$ corresponding to a Grothendieck topology $K$ on $\cal C$ containing $J$. Then $k$ is implicationally open if and only if for any two $J$-closed sieves $S$ and $T$ on an object $c\in {\cal E}$, if for every arrow $f:d\to c$, $f^{\ast}(S)\in J(d)$ implies $f^{\ast}(T)\in J(d)$ then the sieve $\{f:e\to c \textrm{ | } f^{\ast}(S)\subseteq f^{\ast}(T)\}$ is $J$-covering.
\end{proposition}

\begin{proofs}
The thesis follows immediately from the lemma by using the explicit descriptions of the Heyting algebra operations on the subobject lattice $\Omega$ of $[{\cal C}^{\textrm{op}}, \Set]$ obtained in section \ref{sec2}.
\end{proofs}

\begin{theorem}
Let $\cal E$ be a locally small topos and $\cal C$ be a set of objects of $\cal E$ such that the class of objects of $\cal E$ which admit a monomorphism to an object of $\cal C$ forms a separating set for $\cal E$. Then for any geometric morphism $f:{\cal F}\to {\cal E}$, $f^{\ast}$ preserves the Heyting negation (resp. the Heyting implication) operation on subobjects if and only if for every $c\in {\cal C}$, $f^{\ast}:\Sub_{\cal E}(c)\to \Sub_{\cal F}(f^{\ast}(c))$ preserves the Heyting negation (resp. the Heyting implication) operation.
\end{theorem}

\begin{proofs}
Given an epimorphism $p:\coprod_{i\in I}c_{i}\epi e$ in $\cal E$, the pullback functor $p^{\ast}:\Sub_{{\cal E}}(e)\to \Sub_{\cal E}(\coprod_{i\in I}c_{i})\cong \prod_{i\in I}\Sub_{\cal E}(c_{i})$ is logical and conservative, and for any monomorphism $m:b\mono a$ in $\cal E$ the pullback functor $m^{\ast}:\Sub_{\cal E}(a)\to \Sub_{\cal E}(b)$ is logical and essentially surjective (cf. the proof of Theorem \ref{cartesian}). Our thesis thus follows immediately from the fact that the inverse image functors of geometric morphisms preserve pullbacks and arbitrary coproducts.
\end{proofs}

\begin{corollary}
Let $f:L\to L'$ be a morphism of locales, with corresponding geometric morphism $\Sh(f):\Sh(L)\to \Sh(L')$. Then ${\Sh(f)}^{\ast}$ preserves the pseudocomplementation (resp. the Heyting implication) operation on subobjects if and only if $f:L'\to L$ preserves the operation of Heyting negation (resp. of Heyting implication). In particular, for any locale $L$ and nucleus $j$ on $L$, with fixset $L_{j}$, the geometric inclusion $\Sh(L_{j})\hookrightarrow \Sh(L)$ preserves the pseudocomplementation (resp. the Heyting implication) operation on subobjects if and only if $j:L\to L_{j}$ preserves the operation of Heyting negation (resp. of Heyting implication).  
\end{corollary}

This results motivates the following definition.

\begin{definition}
\begin{enumerate}[(a)]
\item A morphism $f:L\to L'$ of locales is said to be \emph{weakly open} (resp. \emph{implicationally open}) if the frame homomorphism $f:L'\to L$ preserves the pseudocomplementation (resp. the Heyting implication) operation;

\item A sublocale of a locale $L$, with corresponding nucleus $j$ on $L$, is said to be \emph{weakly open} (resp. \emph{implicationally open}) if $j:L\to L_{j}$ preserves the pseudocomplementation (resp. the Heyting implication) operation. 
\end{enumerate} 
\end{definition}

\begin{remark}
Let $L$ be a locale and $j$ be a nucleus on $L$. Then
\begin{enumerate}[(i)]
\item The sublocale of $L$ corresponding to $j$ is weakly open if and only if for any $a,b\in L$, $j(a\wedge b)=j(0)$ implies $j(b)=j(b\wedge \neg a)$, where $\neg$ denotes the Heyting pseudocomplementation operation on $L$;  

\item The sublocale of $L$ corresponding to $j$ is implicationally open if and only if for any $a,b, c\in L$, $j(c\wedge a)\leq j(b)$ implies $j(c)\leq j(a\imp b)$, where $\imp$ denotes the operation of Heyting implication on $L$.  
\end{enumerate}
\end{remark}

Recall that a sublocale of a locale $L$ is said to be \emph{dense} if the corresponding nucleus $j$ on $L$ satisfies $j(0)=0$, where $0$ is the bottom element of $L$, equivalently if $j\leq \neg\neg$. Note that any dense sublocale is weakly open; indeed, if $j(a\wedge b)=j(0)=0$ then $a\wedge b\leq j(a\wedge b)=0$ and hence $b\leq \neg a$, which implies $j(b)\leq j(\neg a)$, as required.

\section{Analogues of the double negation and De Morgan topologies}\label{sec4}

Recall from \cite{OC3} that on every topos there exists a local operator $m$, called the De Morgan topology on $\cal E$, satisfying the following properties:

\begin{enumerate}[(i)]
\item $m\leq \neg\neg$;

\item $\sh_{m}({\cal E})$ satisfies De Morgan's law;

\item For every $j\leq \neg\neg$, $j\geq m$ if and only if $\sh_{j}({\cal E})$ satisfies De Morgan's law.
\end{enumerate}

Analogously, one can prove that the Boolean or double negation topology $b=\neg\neg$ on a topos $\cal E$ enjoys the following properties:

\begin{enumerate}[(i)]
\item $b\leq \neg\neg$;

\item $\sh_{b}({\cal E})$ is Boolean;

\item For every $j\leq \neg\neg$, $j\geq b$ if and only if $\sh_{j}({\cal E})$ is Boolean.
\end{enumerate}

In fact, as we shall see in Proposition \ref{relax}, the condition $j\leq \neg\neg$ in property $(iii)$ can be relaxed to the condition that $j$ be weakly open; this remark will pave the way for introducing natural analogues of the double negation and De Morgan topologies for a wide class of intermediate logics.

Recall from \cite{OC7} (Proposition 6.2) that the Boolean or double negation topology $b$ on $\cal E$ is the smallest topology $j$ on $\cal E$ such that all the subobjects in $\cal E$ of the form $m\vee \neg m$ are $j$-dense, equivalently the smallest topology $j$ on $\cal E$ such that the equalizer of the pair of maps $f,g:\Omega\to \Omega$ where $f$ is the composite $\vee \circ <1_{\Omega}, \neg>$ and $g$ is the composite $\top \circ !_{\Omega}$ (where $!_{\Omega}$ is the unique arrow $\Omega\to 1$) is $j$-dense; analogously, the De Morgan topology $m$ on $\cal E$ is the smallest topology $j$ on $\cal E$ such that all the subobjects in $\cal E$ of the form $\neg m\vee \neg \neg m$ are $j$-dense, equivalently the smallest topology $j$ on $\cal E$ such that the equalizer of the pair of maps $f,g:\Omega\to \Omega$ where $f$ is the composite $\vee \circ <\neg, \neg \circ \neg>$ and $g$ is the composite $\top \circ !_{\Omega}$ is $j$-dense. 

\begin{proposition}\label{relax}
Let $\cal E$ be a topos and $k$ be the Boolean (resp. De Morgan) topology on $\cal E$. Then for any weakly open local operator $j$ on $\cal E$, $j\geq k$ if and only if $\sh_{j}({\cal E})$ is Boolean (resp. satisfies De Morgan's law).
\end{proposition}

\begin{proofs}
This follows immediately from the characterization of the Boolean (resp. De Morgan) topology on $\cal E$ as the smallest topology $j$ on $\cal E$ such that all the subobjects in $\cal E$ of the form $m\vee \neg m$ (resp. of the form $\neg m\vee \neg \neg m$) are $j$-dense and the fact that $\cal E$ is Boolean (resp. De Morgan) if and only if all the subobjects in $\cal E$ of the form $m\vee \neg m$ (resp. of the form $\neg m\vee \neg \neg m$) are isomorphisms. Indeed, $a_{j}$ preserves (unions and) pseudocomplementations of subobjects and if $j\geq k$ then any monomorphism which is $k$-dense is also $j$-dense.  
\end{proofs}

Note that the essential point in the proof of Proposition \ref{relax} is the fact that if $j$ is weakly open then the associated sheaf functor $a_{j}$ preserves the `logical structure' of the subobjects involved in the definition of the Boolean (resp. De Morgan) topologies, i.e. $a_{j}(m\vee \neg m)=a_{j}(m)\vee \neg a_{j}(m)$ (resp. $a_{j}(\neg m\vee \neg \neg m)=\neg a_{j}(m)\vee \neg \neg a_{j}(m)$). 

To obtain an analogue of these results for other intermediate logics, we observe that any intermediate logic obtained from intuitionistic logic by adding a single axiom scheme can be associated to a sequent of the form $\top \vdash_{\vec{x}} \top=\phi$, where $\phi(\vec{x})$ is a term in the algebraic theory of Heyting algebras, whose signature consists of two constant symbols $0$ and $1$, one unary operation $\neg$ and three binary operations $\vee$, $\wedge$ and $\imp$. We can define the notion of a topos $\cal E$ satisfying a given intermediate logic $L$ specified by an axiom of the form $\top \vdash_{\vec{x}} \phi_{L}$ as above by requiring that the internal Heyting algebra in $\cal E$ given by its subobject classifier $\Omega$ should satisfy $\phi_{L}$. This condition can be expressed by saying that the equalizer of the pair of maps $f_{L}, g_{L}:{\Omega}^{n}\to \Omega$ (where $n$ is the number of free variables occurring in the term $\phi_{L}(\vec{x})$), where $f_{L}$ is the arrow which represents the interpretation of the term $\phi(\vec{x})$ in the Heyting algebra $\Omega$ and $g_{L}$ is the composite of the unique arrow $!_{\Omega^{n}}:\Omega^{n}\to 1$ with the truth value $\top:1\to \Omega$, is an isomorphism.

In view of the definition of double negation and De Morgan topology, this characterization leads us to defining the notion of $L$-topology on a topos as follows (recall from \cite{El}, specifically pp. 212-214, that the smallest local operator $j$ on a topos for which a given monomorphism is $j$-dense always exists, so the following definition actually makes sense).

\begin{definition}\label{Ltopology}

\begin{enumerate}[(a)]
\item Let $\cal E$ be an elementary topos and $L$ be an intermediate logic presented by an axiom $\top \vdash_{\vec{x}} \top=\phi_{L}$ as above. The \emph{$L$-topology} $j^{\cal E}_{L}$ on $\cal E$ is the smallest local operator $j$ on $\cal E$ such that the equalizer of $f_{L}$ and $g_{L}$ in $\cal E$ is $j$-dense;

\item Let $A$ be a locale and $L$ be an intermediate logic presented by an axiom $\top \vdash_{\vec{x}} \top=\phi_{L}$ as above. The \emph{$L$-sublocale} of $A$ is the sublocale of $A$ corresponding to the subtopos $\sh_{j^{\Sh(A)}_{L}}(\Sh(A))\hookrightarrow \Sh(A)$.   
\end{enumerate} 
\end{definition}

\begin{remark}\label{equiv}
Composing the classifying maps of subobjects in $\cal E$ with $f_{L}$ gives rise to an operation $\xi_{L}$ on subobjects of $\cal E$; in terms of this operation, the condition that, for a given local operator $j$ on $\cal E$, the equalizer of $f_{L}$ and $g_{L}$ should be $j$-dense can be reformulated by saying that all the subobjects arising as the result of applying the operation $\xi_{L}$ to arbitrary subobjects of $\cal E$ are $j$-dense. Note that a subobject in the image $\xi_{L}(m_{1}, \ldots, m_{n})$ of the operation $\xi_{L}$ is a combination of the subobjects $m_{1}, \ldots, m_{n}$ by using the Heyting algebra constants and connectives exactly as they appear in the term $\phi_{L}$ defining the logic $L$.  
\end{remark}

That our general definition actually represents a natural generalization of the double negation and the De Morgan topologies is shown by the following result, which represents a generalization of Proposition \ref{relax} to a wide class of intermediate logics.

\begin{proposition}\label{interm}
Let $\cal E$ be a topos, $L$ be an intermediate logic specified by an axiom $\top \vdash_{\vec{x}} \top=\phi_{L}$ where $\phi_{L}$ is a term obtained by applying the operation $\vee$ to terms only involving the connectives $\wedge$, $\imp$ and $\neg$, and $j^{\cal E}_{L}$ be the $L$-topology on $\cal E$ as in Definition \ref{Ltopology}. Then
\begin{enumerate}[(i)]
\item $j^{\cal E}_{L}\leq \neg\neg$;

\item $\sh_{j^{\cal E}_{L}}({\cal E})$ satisfies $L$;

\item For every implicationally open (resp. weakly open, if the connective $\imp$ does not appear in the term $\phi_{L}$) local operator $j$ on $\cal E$, $j\geq j^{\cal E}_{L}$ if and only if $\sh_{j}({\cal E})$ satisfies $L$.
\end{enumerate}
\end{proposition}

\begin{proofs}
Condition $(i)$ is satisfied by definition of $j^{\cal E}_{L}$ since $\neg\neg$ is an implicationally open operator (cf. Example 4.5.9 \cite{El}) and $\sh_{\neg\neg}({\cal E})$ is Boolean whence it satisfies $L$ ($L$ being by our hypothesis an intermediate logic).

To prove condition $(ii)$ we observe that for any local operator $j$ and any objects $A, B\in \sh_{j}({\cal E})$, the pseudocomplementation $A\imp_{j} B$ in $\sh_{j}({\cal E})$ coincides with the pseudocomplementation $A\imp B$ in $\cal E$. Indeed, the diagram 
\[  
\xymatrix {
\Omega_{j}\times\Omega_{j} \ar[d]^{e_{j}\times e_{j}} \ar[rr]^{\imp_{j}} & & \Omega_{j} \ar[d]^{e_{j}} \\
\Omega\times\Omega \ar[rr]_{\imp} & & \Omega}
\]
commutes. Moreover, if $j$ is dense then the initial object $0$ of $\cal E$ is a $j$-sheaf and hence for any subobject $A\mono E$ in $\sh_{j}({\cal E})$, its pseudocomplementation $\neg_{j}A$ in $\sh_{j}({\cal E})$ coincides with its pseudocomplementation $\neg A$ in $\cal E$. Therefore for any dense local operator $j$ the term $\phi_{L}$ evaluated at any subobject in $\sh_{j}({\cal E})$ (where the connectives $\wedge$, $\vee$, $\imp$ and $\neg$ are respectively interpreted by the Heyting operations of subobjects $\wedge_{j}$, $\vee_{j}$, $\imp_{j}$ and $\neg_{j}$ in $\sh_{j}({\cal E})$) is equal to the result of applying the associated sheaf functor $a_{j}$ to the result of applying $\phi_{L}$ to the same subobjects, but regarded as subobjects in $\cal E$ (by the particular form of the formula $\phi_{L}$, since $a_{j}$ preserves arbitrary unions, that is the interpretation of disjunctions). Therefore, since $j^{\cal E}_{L}$ is dense (by condition $(i)$), by the definition of $L$-topology the sequent $\top \vdash_{\vec{x}} \top=\phi_{L}$ holds in $\sh_{j^{\cal E}_{L}}({\cal E})$ when it is evaluated at any subobject in $\sh_{j^{\cal E}_{L}}({\cal E})$, that is the topos $\sh_{j^{\cal E}_{L}}({\cal E})$ satisfies the logic $L$.       

Condition $(iii)$ follows from the fact that, since $a_{j}$ preserves the `logical structure' of the term $\phi_{L}$, it sends the equalizer of $f_{L}$ and $g_{L}$ in $\cal E$ to an isomorphism (that is, by definition of $L$-topology, $j\geq j^{\cal E}_{L}$) if and only if the topos $\sh_{j}({\cal E})$ satisfies the logic $L$. 
\end{proofs} 

Notice that classical logic, De Morgan logic, G\"{o}del-Dummett logic and Smetanich logic all satisfy the hypotheses of the Proposition \ref{interm}.

The following proposition represents a generalization of Proposition 2.5 \cite{OC3}.

\begin{proposition}
Let $\cal E$ be an elementary topos and $L$ an intermediate logic as in the hypotheses of Proposition \ref{interm}. Then for any implicationally open (resp. weakly open, if the connective $\imp$ does not appear in the term $\phi_{L}$) local operator $j$ on $\cal E$, we have $\sh_{j^{\sh_{j}({\cal E})}_{L}}(\sh_{ j^{\cal E}_{L}}({\cal E}))=\sh_{j^{\cal E}_{L} \vee j}({\cal E})$. 
\end{proposition}

\begin{proofs}
We have to show that the relativization $k$ at $j$ of the local operator $j^{\cal E}_{L} \vee j$ (in the sense of section 6 of \cite{OC6}) satisfies the universal property in the definition of $L$-topology on $\sh_{j}({\cal E})$, that is for any local operator $s$ on $\sh_{j}({\cal E})$, the equalizer $E$ of $f_{L}, g_{L}$ in $\sh_{j}({\cal E})$ is $c_{s}$-dense if and only if $s\geq k$.
Now, given a local operator $s$ on $\sh_{j}({\cal E})$, since the associated sheaf functor $a_{j}:{\cal E}\to \sh_{j}({\cal E})$ preserves the logical structure of the term $\phi_{L}$, $a_{s}$ sends $E$ to an isomorphism if and only if the associated sheaf functor ${\cal E}\to \sh_{s}(\sh_{j}({\cal E}))=\sh_{\tilde{s}}({\cal E})$ corresponding to the local operator $\tilde{s}$ on $\cal E$ whose relativization at $j$ is $s$ sends the equalizer of $f_{L}$ and $g_{L}$ in $\cal E$ to an isomorphism, that is if and only if $\tilde{s}\geq j^{\cal E}_{L}$, equivalently $s\geq k$.
\end{proofs}

One might naturally wonder if there exist explicit formulae for calculating $L$-topologies on Grothendieck toposes. We shall devote the rest of this section to addressing this problem.

Given a Grothendieck topos $\Sh({\cal C}, J)$ and an intermediate logic $L$ corresponding to a term $\phi_{L}$ in the theory of Heyting algebras as above in this section, we can explicitly describe the Grothendieck topology $J_{L}\supseteq J$ on $\cal C$ characterized by the property that the canonical inclusion $\Sh({\cal C}, J_{L})\hookrightarrow \Sh({\cal C}, J)$ can be identified with the subtopos $\sh_{j_{L}}({\Sh({\cal C}, J)})$ of $\Sh({\cal C}, J)$ (where $j_{L}$ is the $L$-topology on $\Sh({\cal C}, J)$), as follows. By definition of $L$-topology, $j_{L}$ is the smallest local operator on $\Sh({\cal C}, J)$ such that the equalizer $E$ of the arrows $f_{L}, g_{L}:{\Omega_{J}}^{n}\to \Omega_{J}$ is dense for it (where $\Omega_{J}$ is the subobject classifier of $\Sh({\cal C}, J)$). Now, if $K$ is the Grothendieck topology on $\cal C$ corresponding to a local operator $j$ on $\Sh({\cal C}, J)$ the condition for $E$ to be $c_{j}$-dense can be expressed by saying that $E$, regarded as a subobject in the presheaf topos $[{\cal C}^{\textrm{op}}, \Set]$, is $c_{K}$-dense, and hence reformulated as the requirement that for any $(S_{1}, \ldots, S_{n})\in {\Omega_{J}}^{n}(c)$, $\{f:d\to c \textrm{ | } f_{L}(f^{\ast}(S_{1}, \ldots, S_{n}))=M_{d}\}\in J(c)$ (where $M_{d}$ is the maximal sieve on $d$). But, since in any topos the pullback functor along any arrow is logical and commutes with any closure operation, $\{f:d\to c \textrm{ | } f_{L}(f^{\ast}(S_{1}, \ldots, S_{n}))=M_{d}\}=\{f:d\to c \textrm{ | } f^{\ast}(f_{L}(S_{1}, \ldots, S_{n}))=M_{d}\}$, that is $\{f:d\to c \textrm{ | } f_{L}(f^{\ast}(S_{1}, \ldots, S_{n}))=M_{d}\}\in J(c)$ if and only if $f_{L}(S_{1}, \ldots, S_{n})\in J(c)$. We can thus conclude that the Grothendieck topology $K$ on ${\cal C}$ corresponding to the $L$-topology on the topos $\Sh({\cal C}, J)$ is the Grothendieck topology on $\cal C$ generated over $J$ by the sieves of the form $f_{L}(S_{1}, \ldots, S_{n})$ (for any $c\in {\cal C}$ and any $(S_{1}, \ldots, S_{n})\in {\Omega_{J}}^{n}(c)$). 

If the formula $\phi_{L}$ satisfies the hypotheses of Proposition \ref{interm} and the topology $J$ on $\cal C$ is dense then one can equivalently describe $K$ as the Grothendieck topology on $\cal C$ generated over $J$ by the sieves of the form $\tilde{f_{L}}(S_{1}, \ldots, S_{n})$ (for any $c\in {\cal C}$ and any $(S_{1}, \ldots, S_{n})\in {\Omega_{J}}^{n}(c)$), where $\tilde{f_{L}}$ is the interpretation of the term $\phi_{L}$ in the presheaf topos $[{\cal C}^{\textrm{op}}, \Set]$ applied the $J$-closed sieves $S_{1}, \ldots , S_{n}$.          

Notice that a standard way of obtaining a dense geometric inclusion from a topos $\Sh({\cal C}, J)$ to a presheaf topos is the following. For any site $({\cal C}, J)$, the full subcategory $\tilde{{\cal C}}$ of $\cal C$ on the objects $c$ such that $\emptyset\in J(c)$ is $J$-dense and hence by Grothendieck's Comparison Lemma we have an equivalence $\Sh({\cal C}, J)\simeq \Sh(\tilde{{\cal C}}, J|_{\tilde{{\cal C}}})$; the induced topology $J|_{\tilde{{\cal C}}}$ on $\tilde{{\cal C}}$ is dense, and hence, if the logic $L$ satisfies the hypotheses of Proposition \ref{interm}, the $L$-topology on $\Sh({\cal C}, J)$ admits a simpler description as the Grothendieck topology on $\tilde{{\cal C}}$ generated over $J|_{\tilde{{\cal C}}}$ by the collection of sieves obtained by applying the term $\phi_{L}$ to $J|_{\tilde{{\cal C}}}$-closed sieves using the interpretation of connectives in the presheaf topos $[\tilde{{\cal C}}^{\textrm{op}}, \Set]$ (rather than in the topos $\Sh({\cal C}, J)$).
        
Let us now apply Proposition \ref{interm} in the context of locales.

\begin{proposition}\label{intermlocales}
Let $A$ be a locale, $L$ be an intermediate logic satisfying the hypotheses of Proposition \ref{interm} above, and $A_{L}$ be the $L$-sublocale of $A$ as in Definition \ref{Ltopology}. Then
\begin{enumerate}[(i)]
\item $A_{L}$ is a dense sublocale of $A$;

\item $A_{L}$ satisfies $L$;

\item For every implicationally open (resp. weakly open, if the connective $\imp$ does not appear in the term $\phi_{L}$) sublocale $A'$ of $A_{L}$, $A'\subseteq A_{L}$ if and only if $A'$ satisfies $L$.
\end{enumerate}
\end{proposition}

Recall that for any frame $A$ and any filter $F$ on $A$, the relation $\simeq$ on $A$ defined by `$a\simeq b$ if and only if $(a\imp b)\wedge (b\imp a)\in F$' is an equivalence relation on $A$, and the quotient set $A\slash F$ is a frame satisfying the universal property that any frame homomorphism $A\to B$ sending every element of $F$ to the top element $1$ of $B$ factors uniquely through the canonical projection $A\to A\slash F$. Thus, under the hypotheses of Proposition \ref{intermlocales}, we can concretely build the $L$-sublocale of a locale $A$ as the quotient of $A$ by the filter generated by the set of elements of the form $\phi_{L}(x_{1}, \ldots, x_{n})$ for $x_{1}, \ldots, x_{n}\in A$. Of course, alternative constructions of $L$-sublocales of given locales, based on different representations of the corresponding toposes, are also possible. For example, one can obtain a different description of the DeMorganization of a locale $A$ (i.e., of the $L$-sublocale of a locale $A$ where $L$ is De Morgan logic, cf. \cite{OC3}) by representing the topos $\Sh(A)=\Sh(A, J_{A})$ (where $J_{A}$ is the canonical topology on $A$) as $\Sh(A^{\ast}, J_{A}|_{A^{\ast}})$, where $A^{\ast}$ is the full subcategory of $A$ on its non-zero elements. Recall from \cite{OC3} that the De Morgan topology on a small category $\cal C$, that is the Grothendieck topology on $\cal C$ corresponding to the $L$-subtopos of $[{\cal C}^{\textrm{op}}, \Set]$ for $L$ equal to De Morgan logic (or equivalently to a dense subtopos $\Sh({\cal C}, J)$ of $[{\cal C}^{\textrm{op}}, \Set]$) is generated by the pullback-stable family of sieves of the form 
\[
\{f:d\rightarrow c \textrm{ | } f^{\ast}(R)=\emptyset \textrm{ or } f^{\ast}(R) \textrm{ is stably non-empty} \}
\]       
where $R$ is a sieve in $\cal C$ on $c$; moreover, if $\cal C$ is a geometric category, one can suppose, without loss of generality, $R$ to be generated by a single arrow. 

From this we easily deduce that the DeMorganization of a locale $A$ can be identified with the surjective frame homomorphism $\xi:A\to A_{m}$ defined as follows: $A_{m}$ is the set of elements $l$ of $A$ with the property that for any elements $r,a \in A$ with $r\leq a$, if $b\leq l$ for every $b\leq a$ such that either $b\wedge r=0$ or for every $c\neq 0$ such that $c\leq b$, $c\wedge r\neq 0$, then $a\leq l$, while $\xi$ is the map sending any $l\in A$ to the smallest element $l'\geq l$ belonging to $A_{m}$. 

Similarly, one can calculate the $L$-sublocale $\xi':A\to A_{GD}$ of a given locale $A$ where $L$ is G\"{o}del-Dummett logic. In general, the Grothendieck topology on a small category $\cal C$ corresponding to the $L$-subtopos of $\Sh({\cal C}, J)$ (where $L$ is G\"{o}del-Dummett logic), is generated by the pullback-stable family of sieves of the form 
\[
\{f:d\to c \textrm{ | } f^{\ast}(R)\subseteq f^{\ast}(S) \textrm{ or } f^{\ast}(S)\subseteq f^{\ast}(R)\}
\] 
for any $J$-closed sieves $R$ and $S$ on an object $c\in {\cal C}$. From this we immediately deduce, similarly to above, that $A_{GD}$ can be described as the set of elements $l$ of $A$ with the property that for any elements $r,s,c\in A$ such that $r,s\leq c$ if $(r\imp s)\wedge c\leq l$ and $(s\imp r)\wedge c\leq l$ then $c\leq l$, while $\xi'$ can be identified with the map sending any element $l\in A$ to the smallest element $l'\geq l$ which belongs to $A_{GD}$.

\end{document}